\documentclass[12pt]{article}

\input{tcilatex}
\begin{document}

\begin{center}
\textbf{Identification des param\`{e}tres d'un mod\`{e}le logistique en
dynamique des populations avec sortie affine}

\ 

Messaoud SOUILAH$^{1}$ \& Imene Sabira SOUALAH$^{2}$

$^{1}$Faculty of Mathematics, USTHB, Algiers

P.O. Box 32 El Alia Bab Ezzouar Algiers 16111 Algeria

msouilah@usthb.dz

$^{2}$High National School in Advenced Technologies

Faculty of Medecine, Biomedical Center, Bordj El Kiffan, Algiers Algeria

sabiraimene.soualah@ensta.edu.dz

02 Octobre 2024
\end{center}

\textbf{ABSTRACT.} We study the parameters identification of a dynamic model
of a population living in a given host environment governed by a logistic
law. We use a statistic Kullback-Leibler type method to derive the algorithm
for estimating the parameters of the model in two levels. The first level
AIG is an offline algorithm and global it is obtained using the critical
points of the reestimation transformation between two parameters. It
estimates the parameters in a global iterative manner starting from a block
of data. The second level ARE is adaptive recursive and is used online. It
constitutes a refinement of the AIG algorithm. The convergence of the AIG
algorithm is an open problem. The convergence of the ARE algorithm is
demonstrated by constructing a new model, a new space and a new probability
law.


\section{Le mod\`{e}le dynamique continu}

Etant donn\'{e}e une population vivant dans un milieu d'accueil $\Omega $
caracteris\'{e}e par deux constantes sp\'{e}ecifiques internes : le taux de
natalit\'{e} $b$ et le taux de mortalit\'{e} $d.$ On notera par $r=b-d$. La
constante $K$ d\'{e}signe la capacit\'{e} d'accueil du milieu, i.e. le
nombre maximal d'individus pouvant vivre dans ce milieu.

La dynamique de cette population, i.e. la loi de variation au cours du temps
est donn\'{e}e par une \'{e}quation d'\'{e}tat et une \'{e}quation
d'observation.

\begin{equation}
\left \{ 
\begin{array}{l}
u_{t}=ru(1-\dfrac{u}{K}) \\ 
y=au+b%
\end{array}%
\right.  \label{md}
\end{equation}

L'\'{e}quation d'\'{e}tat est une loi logistique de la population et l'\'{e}%
quation d'observation est une sortie affine du mod\`{e}le.

\section{Mod\`{e}le dynamique discret}

Soit $h>0$ un pas de discr\'{e}tisation du temps. Le mod\`{e}le dynamique
discret est donn\'{e} par

\begin{equation}
\left \{ 
\begin{array}{l}
u_{n+1}=u_{n}+hru_{n}(1-K^{-1}u_{n}) \\ 
y_{n}=au_{n}+b%
\end{array}%
\right.  \label{mdd}
\end{equation}

\section{Mod\`{e}le dynamique stochastique discret}

Pour compenser diff\'{e}rentes erreurs de mod\'{e}lisation et de discr\'{e}%
tisation, on introduit des bruits blancs additifs dans les deux \'{e}quations

\begin{equation}
\left \{ 
\begin{array}{l}
u_{n+1}=u_{n}+hru_{n}(1-K^{-1}u_{n})+v_{n} \\ 
y_{n}=au_{n}+b+w_{n}%
\end{array}%
\right.  \label{mdsd}
\end{equation}

$v_{n}$ et $w_{n}$ sont des variables al\'{e}atoires Gaussiennes de
variances respectives $\sigma $ et $\tau .$ Les param\`{e}tres du mod\`{e}le 
\`{a} estimer sont les param\`{e}tres \'{e}cologiques $r,K,a,b$ auxquelles
sont ajout\'{e}es les param\`{e}tres math\'{e}matiques $\sigma ,\tau $ dues 
\`{a} la m\'{e}thode stochastique utilis\'{e}e. Le vecteur des param\`{e}%
tres \`{a} estimer est donc

$\lambda =(r,K,a,b,\sigma ,\tau )$

\section{Estimation des param\`{e}tres}

On utilise une m\'{e}thode d'estimation de type HMM (mod\`{e}les de markov
cach\'{e}es). La m\'{e}thode d'estimation est bas\'{e}e sur la
transformation de r\'{e}estimation suivante

\begin{equation}
Q(\lambda ,\lambda ^{\prime })=\sum
\limits_{u_{1}^{N}}L(y_{1}^{N},u_{1}^{N}/\lambda )\ln
L(y_{1}^{N},u_{1}^{N}/\lambda ^{\prime })  \label{kl}
\end{equation}

ou la vraisemblance (Likelihood) $L(y_{1}^{N},u_{1}^{N}/\lambda )$ est d\'{e}%
finie, \`{a} l'aide des \'{e}chantillions des \'{e}tats $u_{1}^{N}$ et
d'observations $y_{1}^{N},$ par

\begin{equation}
L(y_{1}^{N},u_{1}^{N}/\lambda ^{\prime })=P(u_{0})\prod
\limits_{n=1}^{N}P(y_{n}/u_{n},\lambda )P(u_{n}/u_{n-1},\lambda )  \label{vr}
\end{equation}

en fonction des probabilit\'{e}s de transition

\begin{equation}
P(u_{n}/u_{n-1},\lambda )=\dfrac{1}{\sigma \sqrt{2\pi }}\exp (-\dfrac{%
(u_{n}-u_{n-1}-hru_{n-1}(1-K^{-1}u_{n-1}))^{2}}{2\sigma ^{2}})  \label{pob}
\end{equation}

et les probabilit\'{e}s des observations

\begin{equation}
P(y_{n}/u_{n},\lambda )=\dfrac{1}{\tau \sqrt{2\pi }}\exp (-\dfrac{%
(y_{n}-au_{n}-b)^{2}}{2\tau ^{2}})  \label{ptr}
\end{equation}

Cette m\'{e}thode est bas\'{e}e sur l'hypothese d'ind\'{e}pendance
conditionnelle dans les HMM

\begin{equation}
L(y_{1}^{N}/u_{1}^{N},\lambda ^{\prime })=\prod
\limits_{n=1}^{N}P(y_{n}/u_{n},\lambda )
\end{equation}

\section{Algorithme it\'{e}ratif global AIG}

La fonction $-Q(\lambda ,\lambda ^{\prime })$ d\'{e}signe une entropie du
syst\`{e}me (d\'{e}sordre) et est destin\'{e}e \`{a} \^{e}tre minimis\'{e}e,
donc on maximise $Q(\lambda ,\lambda ^{\prime }).$ Pour cela, on utilise le
point fixe du gradient de $\lambda ^{\prime }\mapsto Q(\lambda ,\lambda
^{\prime })$ part rapport a $\lambda ^{\prime }.$ On a

\begin{eqnarray*}
\ln L(y_{1}^{N},u_{1}^{N}/\lambda ^{\prime }) &=&\ln P(u_{0})+ \\
&&\sum \limits_{n=1}^{N}\left( 
\begin{array}{c}
-\ln (\tau ^{\prime }\sqrt{2\pi })-\ln (\sigma ^{\prime }\sqrt{2\pi })-%
\dfrac{(y_{n}-a^{\prime }u_{n}-b^{\prime })^{2}}{2\tau ^{\prime 2}} \\ 
-\dfrac{(u_{n}-u_{n-1}-hr^{\prime }u_{n-1}(1-K^{\prime -1}u_{n-1}))^{2}}{%
2\sigma ^{\prime 2}}%
\end{array}%
\right)
\end{eqnarray*}

On obtient pour les param\`{e}tres $\sigma ,\tau :$

\begin{eqnarray*}
\dfrac{\partial }{\partial \tau ^{\prime }}Q(\lambda ,\lambda ^{\prime })
&=&\sum \limits_{n=1}^{N}\sum \limits_{u_{n}}L(y_{1}^{N},u_{n}/\lambda
^{\prime })\left( -\dfrac{1}{\tau ^{\prime }}+\dfrac{(y_{n}-a^{\prime
}u_{n}-b^{\prime })^{2}}{\tau ^{\prime 3}}\right) \\
\dfrac{\partial }{\partial \sigma ^{\prime }}Q(\lambda ,\lambda ^{\prime })
&=&\sum \limits_{n=1}^{N}\sum \limits_{u_{n}}\sum
\limits_{u_{n-1}}L(y_{1}^{N},u_{n},u_{n-1}/\lambda ^{\prime }) \\
&&\left( -\dfrac{1}{\sigma ^{\prime }}+\dfrac{(u_{n}-u_{n-1}-hr^{\prime
}u_{n-1}(1-K^{\prime -1}u_{n-1}))^{2}}{\sigma ^{\prime 3}}\right)
\end{eqnarray*}

Les points critiques de $Q$ en $\lambda ^{\prime }$ donnent

\begin{eqnarray*}
\tau ^{\prime 2} &=&\dfrac{\sum \limits_{n=1}^{N}\sum
\limits_{u_{n}}L(y_{1}^{N},u_{n}/\lambda ^{\prime })(y_{n}-a^{\prime
}u_{n}-b^{\prime })^{2}}{N.L(y_{1}^{N}/\lambda )} \\
\sigma ^{\prime 2} &=&\dfrac{\sum \limits_{n=1}^{N}\sum \limits_{u_{n}}\sum
\limits_{u_{n-1}}L(y_{1}^{N},u_{n},u_{n-1}/\lambda ^{\prime
})(u_{n}-u_{n-1}-hr^{\prime }u_{n-1}(1-K^{\prime -1}u_{n-1}))^{2}}{%
N.L(y_{1}^{N}/\lambda )}
\end{eqnarray*}

Pour les param\`{e}tres $r,K,$ l'algorithme global est plus compliqu\'{e}

\begin{eqnarray*}
\dfrac{\partial }{\partial r^{\prime }}Q(\lambda ,\lambda ^{\prime }) &=&-%
\dfrac{1}{\sigma ^{\prime 2}}\sum \limits_{n=1}^{N}\sum \limits_{u_{n}}\sum
\limits_{u_{n-1}}L(y_{1}^{N},u_{n},u_{n-1}/\lambda ^{\prime
})hu_{n-1}(1-K^{\prime -1}u_{n-1}) \\
&&\left( u_{n}-u_{n-1}-hr^{\prime }u_{n-1}(1-K^{\prime -1}u_{n-1})\right) \\
\dfrac{\partial }{\partial K^{\prime }}Q(\lambda ,\lambda ^{\prime }) &=&-%
\dfrac{1}{\sigma ^{\prime 2}}\sum \limits_{n=1}^{N}\sum \limits_{u_{n}}\sum
\limits_{u_{n-1}}L(y_{1}^{N},u_{n},u_{n-1}/\lambda ^{\prime }) \\
&&(u_{n}-(1+hr^{\prime })u_{n-1}+hr^{\prime }K^{\prime -1}u_{n-1}^{2})^{2}
\end{eqnarray*}

Les points critiques de $Q$ donnent

\begin{eqnarray*}
r^{\prime } &=&\dfrac{\sum \limits_{n=1}^{N}\sum \limits_{u_{n}}\sum
\limits_{u_{n-1}}L(y_{1}^{N},u_{n},u_{n-1}/\lambda ^{\prime })u_{n-1}\left(
u_{n}-u_{n-1}\right) (1-K^{\prime -1}u_{n-1})}{h\sum \limits_{n=1}^{N}\sum
\limits_{u_{n}}\sum \limits_{u_{n-1}}L(y_{1}^{N},u_{n},u_{n-1}/\lambda
^{\prime })u_{n-1}^{2}(1-K^{\prime -1}u_{n-1})^{2}} \\
K^{\prime -1} &=&-\dfrac{\sum \limits_{n=1}^{N}\sum \limits_{u_{n}}\sum
\limits_{u_{n-1}}L(y_{1}^{N},u_{n},u_{n-1}/\lambda ^{\prime
})(u_{n}-u_{n-1}-hr^{\prime }u_{n-1})}{\sum \limits_{n=1}^{N}\sum
\limits_{u_{n}}\sum \limits_{u_{n-1}}L(y_{1}^{N},u_{n},u_{n-1}/\lambda
^{\prime })hr^{\prime }u_{n-1}^{2}}
\end{eqnarray*}

Pour les parametres $a,b,$ l'algorithme global est de complexite moindre

\begin{eqnarray*}
\dfrac{\partial }{\partial a^{\prime }}Q(\lambda ,\lambda ^{\prime }) &=&%
\dfrac{1}{\tau ^{\prime 2}}\sum \limits_{n=1}^{N}\sum
\limits_{u_{n}}L(y_{1}^{N},u_{n}/\lambda ^{\prime })(y_{n}-a^{\prime
}u_{n}-b^{\prime })u_{n} \\
\dfrac{\partial }{\partial b^{\prime }}Q(\lambda ,\lambda ^{\prime }) &=&%
\dfrac{1}{\tau ^{\prime 2}}\sum \limits_{n=1}^{N}\sum
\limits_{u_{n}}L(y_{1}^{N},u_{n}/\lambda ^{\prime })(y_{n}-a^{\prime
}u_{n}-b^{\prime })
\end{eqnarray*}

Les points critiques de $Q$ donnent

\begin{eqnarray*}
a^{\prime } &=&\dfrac{\sum \limits_{n=1}^{N}\sum
\limits_{u_{n}}L(y_{1}^{N},u_{n}/\lambda ^{\prime })(y_{n}-b^{\prime })u_{n}%
}{\sum \limits_{n=1}^{N}\sum \limits_{u_{n}}L(y_{1}^{N},u_{n}/\lambda
^{\prime })u_{n}^{2}} \\
b^{\prime } &=&\dfrac{\sum \limits_{n=1}^{N}\sum
\limits_{u_{n}}L(y_{1}^{N},u_{n}/\lambda ^{\prime })(y_{n}-a^{\prime }u_{n})%
}{N.L(y_{1}^{N}/\lambda ^{\prime })}
\end{eqnarray*}

L'algorithme global pour les six param\`{e}tres est r\'{e}sum\'{e} par

\begin{eqnarray}
r^{\prime } &=&\dfrac{\sum \limits_{n=1}^{N}\sum \limits_{u_{n}}\sum
\limits_{u_{n-1}}L(y_{1}^{N},u_{n},u_{n-1}/\lambda ^{\prime })u_{n-1}\left(
u_{n}-u_{n-1}\right) (1-K^{\prime -1}u_{n-1})}{h\sum \limits_{n=1}^{N}\sum
\limits_{u_{n}}\sum \limits_{u_{n-1}}L(y_{1}^{N},u_{n},u_{n-1}/\lambda
^{\prime })u_{n-1}^{2}(1-K^{\prime -1}u_{n-1})^{2}}  \label{AIG} \\
K^{\prime -1} &=&-\dfrac{\sum \limits_{n=1}^{N}\sum \limits_{u_{n}}\sum
\limits_{u_{n-1}}L(y_{1}^{N},u_{n},u_{n-1}/\lambda ^{\prime
})(u_{n}-u_{n-1}-hr^{\prime }u_{n-1})}{\sum \limits_{n=1}^{N}\sum
\limits_{u_{n}}\sum \limits_{u_{n-1}}L(y_{1}^{N},u_{n},u_{n-1}/\lambda
^{\prime })hr^{\prime }u_{n-1}^{2}}  \nonumber \\
a^{\prime } &=&\dfrac{\sum \limits_{n=1}^{N}\sum
\limits_{u_{n}}L(y_{1}^{N},u_{n}/\lambda ^{\prime })(y_{n}-b^{\prime })u_{n}%
}{\sum \limits_{n=1}^{N}\sum \limits_{u_{n}}L(y_{1}^{N},u_{n}/\lambda
^{\prime })u_{n}^{2}}  \nonumber \\
b^{\prime } &=&\dfrac{\sum \limits_{n=1}^{N}\sum
\limits_{u_{n}}L(y_{1}^{N},u_{n}/\lambda ^{\prime })(y_{n}-a^{\prime }u_{n})%
}{N.L(y_{1}^{N}/\lambda ^{\prime })}  \nonumber \\
\sigma ^{\prime 2} &=&\dfrac{\sum \limits_{n=1}^{N}\sum \limits_{u_{n}}\sum
\limits_{u_{n-1}}L(y_{1}^{N},u_{n},u_{n-1}/\lambda ^{\prime
})(u_{n}-u_{n-1}-hr^{\prime }u_{n-1}(1-K^{\prime -1}u_{n-1}))^{2}}{%
N.L(y_{1}^{N}/\lambda )}  \nonumber \\
\tau ^{\prime 2} &=&\dfrac{\sum \limits_{n=1}^{N}\sum
\limits_{u_{n}}L(y_{1}^{N},u_{n}/\lambda ^{\prime })(y_{n}-a^{\prime
}u_{n}-b^{\prime })^{2}}{N.L(y_{1}^{N}/\lambda )}  \nonumber
\end{eqnarray}

\section{Algorithme r\'{e}cursif}

Ecrivons l'un des six param\`{e}tres dans l'algorithme global AIG
formellement sous la forme $\lambda ^{\prime }=\frac{1}{N}\sum
\limits_{n=1}^{N}S_{n}$

\begin{eqnarray*}
\lambda _{N} &=&\dfrac{1}{N}\sum \limits_{n=1}^{N}S_{n} \\
&=&\dfrac{1}{N}\left( S_{N}+\sum \limits_{n=1}^{N-1}S_{n}\right) \\
&=&\dfrac{1}{N}\left( S_{N}+(N-1)\lambda _{N-1}\right) \\
&=&\left( 1-\dfrac{1}{N}\right) \lambda _{N-1}+\dfrac{1}{N}S_{N}
\end{eqnarray*}

L'algorithme r\'{e}cursif prend donc de la forme

\[
\lambda _{n}=\left( 1-\dfrac{1}{n}\right) \lambda _{n-1}+\dfrac{1}{n}S_{n} 
\]

Il constitue une sorte de combinaison convexe entre la valeur pr\'{e}c\'{e}%
dente $n-1$ et l'histhorique des valeurs de $1$ \`{a} $n.$ Cet algorithme 
\'{e}crit sous cette forme est valable pour les param\`{e}tres $b,\sigma
,\tau .$

Tandis que pour les param\`{e}tres $r,K,a$ il est plus compliqu\'{e} mais va 
\^{e}tre \'{e}crit sous la m\^{e}me forme en rempla\c{c}ant $\dfrac{1}{n}$
par une suite $\varepsilon _{n}$ tendant vers zero dont la s\'{e}rie $\sum
\limits_{n\geq 0}\varepsilon _{n}$ diverge pour des raisons th\'{e}oriques li%
\'{e}s \`{a} la th\'{e}orie des martingales en th\'{e}orie d'estimation
statistique. Par exemple $\varepsilon _{n}=\dfrac{1}{n(\ln n)^{k}},\ k\geq
1. $

\begin{eqnarray}
r_{n} &=&\left( 1-\varepsilon _{n}\right) r_{n-1}+  \label{ARE} \\
&&\varepsilon _{n}\dfrac{\sum \limits_{u_{n}}\sum
\limits_{u_{n-1}}L(y_{1}^{n},u_{n},u_{n-1}/\lambda _{n-1})u_{n-1}\left(
u_{n}-u_{n-1}\right) (1-K_{n-1}^{-1}u_{n-1})}{h\sum \limits_{u_{n}}\sum
\limits_{u_{n-1}}L(y_{1}^{n},u_{n},u_{n-1}/\lambda
_{n-1})u_{n-1}^{2}(1-K_{n-1}^{-1}u_{n-1})^{2}}  \nonumber \\
K_{n}^{-1} &=&\left( 1-\varepsilon _{n}\right) K_{n-1}-  \nonumber \\
&&\varepsilon _{n}\dfrac{\sum \limits_{u_{n}}\sum
\limits_{u_{n-1}}L(y_{1}^{n},u_{n},u_{n-1}/\lambda
_{n-1})(u_{n}-u_{n-1}-hr_{n-1}u_{n-1})}{h\sum \limits_{u_{n}}\sum
\limits_{u_{n-1}}L(y_{1}^{n},u_{n},u_{n-1}/\lambda _{n-1})r_{n-1}u_{n-1}^{2}}
\nonumber \\
a_{n} &=&\left( 1-\varepsilon _{n}\right) a_{n-1}+  \nonumber \\
&&\varepsilon _{n}\dfrac{\sum \limits_{u_{n}}L(y_{1}^{n},u_{n}/\lambda
_{n-1})(y_{n}-b_{n-1})u_{n}}{\sum \limits_{u_{n}}L(y_{1}^{n},u_{n}/\lambda
_{n-1})u_{n}^{2}}  \nonumber \\
b_{n} &=&\left( 1-\varepsilon _{n}\right) b_{n-1}+  \nonumber \\
&&\varepsilon _{n}\dfrac{\sum \limits_{u_{n}}L(y_{1}^{n},u_{n}/\lambda
_{n-1})(y_{n}-a_{n-1}u_{n})}{L(y_{1}^{n}/\lambda _{n-1})}  \nonumber \\
\sigma _{n}^{2} &=&\left( 1-\varepsilon _{n}\right) \sigma _{n-1}^{2}+ 
\nonumber \\
&&\varepsilon _{n}\dfrac{\sum \limits_{u_{n}}\sum
\limits_{u_{n-1}}L(y_{1}^{n},u_{n},u_{n-1}/\lambda
_{n-1})(u_{n}-u_{n-1}-hr_{n-1}u_{n-1}(1-K_{n-1}^{-1}u_{n-1}))^{2}}{%
L(y_{1}^{n}/\lambda _{n-1})}  \nonumber \\
\tau _{n}^{2} &=&\left( 1-\varepsilon _{n}\right) \tau _{n-1}^{2}+  \nonumber
\\
&&\varepsilon _{n}\dfrac{\sum \limits_{u_{n}}L(y_{1}^{n},u_{n}/\lambda
_{n-1})(y_{n}-a_{n-1}u_{n}-b_{n-1})^{2}}{L(y_{1}^{n}/\lambda _{n-1})} 
\nonumber
\end{eqnarray}

\section{Convergence de l'algorithme r\'{e}cursif AER}

La convergence de l'algorithme AIG (\ref{AIG}) est un probl\`{e}me ouvert m%
\^{e}me pour des mod\`{e}les simples. Dans \cite{20},\cite{23} nous avons
conjectur\'{e} qu'il s'agit du th\'{e}or\`{e}me ergodique de Birkhoff dans
le cas g\'{e}n\'{e}ral, mais les non lin\'{e}arit\'{e}s des mod\`{e}les
rendent la d\'{e}monstration techniquement impossible voire impossible. La
convergence de l'algorithme ARE (\ref{ARE}) est d\'{e}montr\'{e}e ici pour
quelques param\`{e}tres parmi les six de $\lambda $ en construisant un
nouveau mod\`{e}le, un nouvel espace et une nouvelle loi de probabilit\'{e}
qui assure la convergence presque s\^{u}re des estimateurs (\ref{ARE}) vers
les vraies valeurs des param\`{e}tres $\lambda .$

Soit $\lambda _{n}$ un param\`{e}tre estim\'{e} du mod\`{e}le non lin\'{e}%
aire par l'algorithme ARE \`{a} l'instant $n$. Cette convergence de
l'algorithme ARE dans les mod\`{e}les non lin\'{e}aires par rapport \`{a} l'%
\'{e}tat $u_{n}$ mais lin\'{e}aires par rapport aux param\`{e}tres $\lambda $
peut \^{e}tre prouv\'{e}e en suivant les \'{e}tapes suivantes \cite{13},\cite%
{14},\cite{15},\cite{20},\cite{21},\cite{22}\cite{23}.

\begin{itemize}
\item \textbf{Etape 1:} Soit $\Omega =\{u_{1}^{N},y_{1}^{N}\},\ N\geq 1,\ 
\mathcal{F}$ une $\sigma $-alg\`{e}bre contenant $\mathcal{F}_{n}:\mathcal{=}%
\sigma (y_{0},\ldots ,y_{n})$ pour tout $n\geq 1$, soit $(\Omega ,\mathcal{F}%
,P)$ un espace de probabilit\'{e} o\`{u} $P$ correspondant aux bruits blancs
additifs et $(\Omega ,\mathcal{F},\widehat{P})$ l'espace de probabilit\'{e} o%
\`{u} $\widehat{P}$ satisfait la condition $\widehat{E}(\widehat{P}%
(y_{n}/u_{n}))=1.$ Alors la suite de variables al\'{e}atoires $(\lambda
_{n}) $ est une martingale sous $\widehat{P}$, i.e. $\widehat{E}(\lambda
_{n}/\mathcal{F}_{n-1}))=\lambda _{n-1}$.

\item \textbf{Etape 2:} Sous les hypoth\`{e}ses pr\'{e}c\'{e}dentes, la
suite $(\lambda _{n})$ est de carr\'{e} int\'{e}grable sous $\widehat{P}$,
i.e. $\widehat{E}(\lambda _{n}^{2})<+\infty $.

\item \textbf{Etape 3:} Les param\`{e}tres estim\'{e}s r\'{e}cursivement par
ARE convergent presque s\^{u}rement sous $\widehat{P}$.

\item \textbf{Etape 4:} Les lois de probabilit\'{e} $P$ relative au mod\`{e}%
le initial et $\widehat{P}$ relative au nouveau mod\`{e}le sont absolument
contigues.

\item \textbf{Etape 5:} Les param\`{e}tres estim\'{e}s r\'{e}cursivement par
ARE convergent presque s\^{u}rement sous $P$.
\end{itemize}

Ces \'{e}tapes ne peuvent pas \^{e}tre appliqu\'{e}es pour toutes les non lin%
\'{e}arit\'{e}s rencontr\'{e}es dans les mod\`{e}les \`{a} cause de la lin%
\'{e}arit\'{e} de l'esp\'{e}rance math\'{e}matique qu'exige la m\'{e}thode.
C'est pour cette raison que nous avons conjectur\'{e} dans \cite{20} que
l'algorithme AIG dont ARE est un raffinement est bas\'{e} sur la th\'{e}orie
ergodique pour prouver cette convergence d'une mani\`{e}re g\'{e}n\'{e}rale.
Mais les d\'{e}tails de cette conjecture exigent un effort math\'{e}matique
consid\'{e}rable.


\newpage

\begin{thebibliography}{99}
\bibitem{10} R.J. Elliot, L. Aggoun, J.B. Moore, Hidden Markov Models.
Estimation and Control, Springer, 1995;

\bibitem{11} G.C. Goodwin, K.S. Sin, Adaptive Filtering Prediction and
Control, Prentice-Hall, Englewood Cliffs, NJ, 1984.

\bibitem{13} A. Khoukhi, T. Aliziane, M. Souilah, Un Algorithme MultiNiveau
d'Identification d'un Canal en Communication Nume%
\'{}%
rique, JESA 36 (4) (2002) 519--537.

\bibitem{14} M. Souilah, A. Khoukhi, T. Aliziane, A new multi-level
algorithm for identification and stochastic adaptive control of industrial
manipulators, Eng. Simulation 26 (4) (2004) 83--98.

\bibitem{15} M. Souilah, A new strategy for identification and control of
mobile robots, Eng. Simulation 28 (3) (2006) 35--48.

\bibitem{16} M. Smordinsky, Ergodic Theory, Entropy, Lecture Notes in
Mathematics, vol. 214, Springer, Berlin, 1971.

\bibitem{17} P. Walters, Lectures on Ergodic Theory, Lecture Notes in
Mathematics, vol. 458, Springer, Berlin, 1975.

\bibitem{18} L.A. Liporace, Maximum likelihood estimation for multivariate
observations of Markov sources, IEEE Trans. Inf. Theory IT-28 (5) (1982)
729--734.

\bibitem{19} R. Zweimu%
\"{}
ller, Hopf's Ratio Ergodic Theorem by Inducing. Communication of R.
Rudnicki. hwww.us.edu.pli.

\bibitem{20} M. Souilah, On some evolution and identification problems in
nonlinear dynamic systems, Doctorate Thesis, Faculty of Mathematics, USTHB
Algiers, 2006.

\bibitem{21} M. Souilah, Parameters Identification in Nonlinear Dynamical
Systems, ODE-SDE School, Faculty of Mathematics, USTHB, 11--16 November,
2006.

\bibitem{22} M. Souilah, Parameters Identification in Nonlinear Dynamic
Systems, TAM-TAM'07, Faculty of Mathematics, USTHB, Tipaza Algiers, 16--18
April 2007, pp. 472--480

\bibitem{23} M. Souilah, A new nonlinear filter for parameters
identification in dynamic systems and application to a transmission channel,
Signal Processing 88, pp. 349--357, 2008.
\end{thebibliography}
\end{document}